\documentclass[12pt,a4paper,twoside]{article}%{amsart}
\usepackage[latin1]{inputenc}
\usepackage[T1]{fontenc}
\usepackage[brazil,english]{babel}
\usepackage{amsfonts}
\usepackage{amssymb}
\usepackage{amsmath}
\usepackage{amsthm}
\usepackage[usenames,dvipsnames]{color}
\usepackage{enumitem}
\usepackage[normalem]{ulem}%provides \uline \uwave \sout
\usepackage{setspace}

\numberwithin{equation}{section}
\theoremstyle{plain} %boldface title, italicized body.
\newtheorem{theorem}{Theorem}[section]
\newtheorem*{theorem*}{Theorem}
\newtheorem{corollary}[theorem]{Corollary}
\newtheorem{lemma}[theorem]{Lemma}
\newtheorem*{lemma*}{Lemma}
\newtheorem{prop}[theorem]{Proposition}
\newtheorem{proposition}[theorem]{Proposition}
\theoremstyle{definition} %boldface title, romand body.

\newtheorem{remark}[theorem]{Remark}
\theoremstyle{remark} % italicized title, romman body.

\newcommand{\eor}{\hfill$\triangleleft$}

\makeatletter
\newcommand{\aslabel}[1]{#1\def\@currentlabel{#1}}
\newcommand{\nowlabel}[1]{\def\@currentlabel{#1}}
\makeatother

\newcommand{\R}{{\mathbb{R}}}

\newcommand{\Cont}{{\mathcal{C}}}

\newcommand{\tref}[1]{\ref{#1}}
\newcommand{\pref}[1]{(\ref{#1})}

\newcommand{\pt}[1]{\left({#1}\right)}
\newcommand{\pq}[1]{\left[{#1}\right]}
\newcommand{\pg}[1]{\left\{{#1}\right\}}
\newcommand{\n}[1]{\left\|{#1}\right\|}
\newcommand{\m}[1]{\left|{#1}\right|}

\usepackage{graphicx}

\definecolor{collnk}{rgb}{.2,0.2,.6}
\definecolor{colcit}{rgb}{.1,0.5,.2}
\usepackage[colorlinks=true]{hyperref} \hypersetup{urlcolor=blue, citecolor=colcit,linkcolor=collnk}

\renewcommand{\l}{\ell}
\newcommand{\s}{\sigma}
\newcommand{\p}{\varpi}
\newcommand{\q}{\ell}
\newcommand{\llambda}{\gamma}
\newcommand{\ttau}{\alpha}
\newcommand{\lt}{\widetilde \l}

\newcommand{\W}{W_0^{1,p}}

\begin{document}

\title {Sobolev versus H\"older local minimizers in degenerate Kirchhoff type problems} 

\author{ \sc{Leonelo Iturriaga}\,$^a$\ \  and\ \ {\sc Eugenio Massa}\,$^{b,}$\footnote{Corresponding author: tel ++55 16 33736635, fax ++55 16 33739650.}
\\[0.4cm]
$^a$ {\small Departamento de Matem\'atica, Universidad T\'ecnica Federico Santa Mar\'ia}
\\[-0.2cm] {\small Avenida Espa\~{n}a 1680, Casilla 110-V, Valpara\'{\i}so, Chile }
\\[-0.2cm] {\small email: leonelo.iturriaga@usm.cl}
\\[0.2cm]
$^b$ \small Departamento de Matem\'atica,
\\[-0.2cm]\small Instituto de Ci\^encias Matem\'aticas e de Computa\c c\~ao, Universidade de S\~ao Paulo,
\\[-0.2cm]\small Campus de S\~ao Carlos, Caixa Postal 668, 13560-970, S\~ao Carlos SP, Brazil.
\\[-0.2cm]\small e-mail: eug.massa@gmail.com 
}

%\date{May 17, 2019 - PREPRINT}

%\author{ PRE-PREPRINT - Friends only \\ Leonelo Iturriaga and   Eugenio Massa.}
\maketitle

\begin{abstract}
In this paper we 
study the geometry of certain functionals associated to quasilinear elliptic boundary value problems with a degenerate nonlocal term of Kirchhoff type.

Due to the degeneration of the nonlocal term it is not possible to directly use classical results such as uniform a-priori estimates and  ``Sobolev versus H\"older local minimizers'' type of results.
We prove that results similar to these  hold true or not, depending on how degenerate the problem is. 

We apply our findings  in order to show existence and multiplicity of solutions for the associated quasilinear equations, considering several different interactions between the nonlocal term and the nonlinearity. 

\medskip

\noindent  {\bf  Mathematical Subject Classification MSC2010:} 
35J20   %	 	Variational methods for second-order elliptic equations
(35J70)   %	Degenerate elliptic equations

\noindent {\bf Key words and phrases:} nonlocal elliptic problems, Kirchhoff equation, degenerate problems, variational methods, a-priori estimates, critical point theory, counterexamples, $p$-Laplacian, local minimizers, existence and multiplicity.

\end{abstract}
%\tableofcontents

\section{Introduction}

The main purpose of this paper is to study the geometry of certain functionals that arise in the study of elliptic equations with a nonlocal term of Kirchhoff type.

Consider for instance the equation
\begin{equation}\label{eq_main_W_f}
\left\{
\begin{array}{ll}
-M(\n{u}_W^{p})\Delta_p u= f(x,u) & \mbox{ in }\ \Omega\,,
\\
u=0 & \mbox{ on }\ \partial \Omega\,,
\end{array}
\right.
\end{equation}
where $\Omega\subset\mathbb{R}^{N}$ is a bounded and smooth domain,  $M$ is some function, 
$\n{\cdot}_W$ is the norm in $W^{1,p}_0(\Omega)$, $p>1$, and $f$ is a suitable nonlinearity.
This problem can be studied variationally, by considering the associated functional 
\begin{equation} \label{eq_funct_gen}
J(u)=\frac1p\widehat M (\n{u}^p_W)-\int_\Omega F(x,u)\,,\qquad u\in W_0^{1,p}(\Omega)\,,
\end{equation}
 where $\widehat M(t)=\int_0^tM(s)\,ds$  and $F(x,v)=\int_0^vf(x,s)\,ds$.

One interesting question that arises in this study is related to the well known results by Brézis  and Nirenberg \cite{H1versusC1}, later extended to the quasilinear case in \cite{ GPM_sob_vs_hold,GuoZhang_sob_vs_hold,BrItUb_multplap}, which  state that a minimum in the $\Cont^1$ topology for a functional in a suitable form, in particular whose  principal term is $\frac1p \n{u}^p_{W}$,  is also a minimum in the $W^{1,p}_0$ topology.
When the principal term of the functional takes the form in \pref{eq_funct_gen}, that is, when considering non local problems as \pref{eq_main_W_f}, these  results may not apply any more. 

In \cite{FanC1W1}, a generalization  of this kind of results to certain Kirchhoff type problems is obtained for the $p(x)$-Laplacian, however, the function $M$ is assumed to be bounded away from zero in compact sets, which rules out the possibility of some degeneracy which might change the picture.  Our aim is to avoid such    nondegeneracy condition, in order to be able to treat a wider class of nonlocal terms $M$. We will focus on the case where the degeneracy happens at the origin, that is, when $M$ is continuous with $M(0)=0$.  As a model one can take $M$ to be a power (see equation \pref{eq_rpq_p}). 

Our main theoretical results are contained in the theorems \tref{th_nomin_gen} and \tref{th_simin}. They show that in some cases, namely when the degeneration of $M$ at the origin is  stronger, no result similar to those in \cite{H1versusC1,GPM_sob_vs_hold,GuoZhang_sob_vs_hold,BrItUb_multplap} can hold true for $J$.  In fact, we provide a situation where the origin is a local minimum when compared to functions that are small in $L^\infty$ or in $\Cont^1$, but  there exists a sequence of functions whose Sobolev norm goes to zero while $J$  stays below its level at  the   origin.   On the other hand, for less degenerate $M$, the above cannot happen, and it is instead possible to prove  that a minimum with respect to the  $L^\infty$ norm is eventually a minimum with respect to the  Sobolev norm. Under further conditions, in Theorem \tref{th_simin_C1}, we also  obtain the full classical result that even a minimum with respect to the  $\Cont^1$ norm is  a minimum with respect to the  Sobolev norm. 

The need for the above distinction between minima with respect to the  $L^\infty$ or  $\Cont^1$ norm is due to the observation that, when dealing with the Kirchhoff operator with degenerate nonlocal term,  regularity estimates like those in \cite{Lieberman} might not hold uniformly, as we show in an example in Section \tref{sec_CE_Lieb} (see Proposition \tref{prop_priori}).  This means that a crucial step of the classical proof of the ``$W^{1,p}$ versus $\Cont^1$'' type results, that is, when one bootstraps from $L^\infty$ estimates to $\Cont^{1,\alpha}$ estimates  in order to apply  Ascoli-Arzela theorem, will be more delicate and not always possible in this context (see Theorem \tref{th_simin_C1}). 
However, in many application, the full strength of the ``$W^{1,p}$ versus $\Cont^1$'' result is not required, and the formulation in terms of the  $L^\infty$ norm is satisfactory: see for example  \cite{Perer_concav,MIO_ODA}. % and the applications in this paper.

As an application, in Section \tref{sec_appl}, we exploit the two results mentioned above
to show existence and multiplicity of solutions for certain nonlocal equations like \pref{eq_main_W_f}. 
It will turn out that the geometry of the associated functional changes  when the degeneracy of the nonlocal term is such that either Theorem \tref{th_simin} or Theorem \tref{th_nomin_gen} holds true, producing different results in terms of existence and multiplicity of the solutions, even if the interaction with the nonlinearity is unchanged: see in particular Theorem \tref{th_coe_mod} and Theorem \tref{th_notpure}. See also the remarks in Section \tref{sec_commap}.

As model examples for the functional \pref{eq_funct_gen} we will consider 
\newcommand{\J}{\mathcal{J}}
\begin{equation}\label{eq_funct_p}
\J(u)=\frac1r \n{u}^r_W+\frac1q\n{u}_q^q-\frac\lambda\p \n{u}_\p^\p\,,
\end{equation}
\begin{equation}\label{eq_funct_p+}
\J^+(u)=\frac1r \n{u}^r_W+\frac1q\n{u^+}_q^q-\frac\lambda\p \n{u^+}_\p^\p\,,
\end{equation}
whose critical points  correspond to solutions of the nonlocal equation 
\begin{equation}\label{eq_rpq_p}
\left\{
\begin{array}{ll}
-\n{u}_W^{r-p}\Delta_p u= -|u|^{q-2}u+\lambda |u|^{\p-2}u & \mbox{ in }\ \Omega\,,
\\
u=0 & \mbox{ on }\ \partial \Omega\,.
\end{array}
\right.
\end{equation}
This problem corresponds to taking, in \pref{eq_funct_gen},
\begin{equation}\label{eq_MF_mod}\begin{cases}
\widehat M(s^p)=\frac pr s^r\,,\\
F(x,v)=\frac\lambda\p |v|^\p-\frac1q |v|^q\quad \text{or} \quad  F(x,v)=\frac\lambda\p (v^+)^\p-\frac1q (v^+)^q\,.
\end{cases}\end{equation}
For $1<q<\p<\frac{pN}{N-p}$ and $r>0$, the origin is a local minimum with respect to the $L^\infty$ norm for the functionals $\J$ and $\J^+$, but a consequence of   Theorem~\tref{th_nomin_gen} and Theorem~\tref{th_simin} will be that such minimum is also a minimum in $\W$ if $r<\frac{pN}{N-p}$, but it is not if $r>\frac{pN}{N-p}$.

\subsection{Some background }

The main feature of problem \pref{eq_main_W_f} is the presence of the term  $M(\n{u}_W^p)$, which is said to be nonlocal, since it depends not only on the point in $\Omega$  where the equation is evaluated, but on the norm of the whole solution. Such problems are  usually called of Kirchhoff type, as they are generalizations of the (stationary) Kirchhoff equation, originally proposed in \cite{Kirch1883} as an improvement of the vibrating string equation, in order to take into account the variation  in the  tension of the string due to the variation of  its  length with respect to the unstrained position. In the  Kirchhoff case, the proposed function  $M$ takes the form  $M(t)=a+bt$ with $a,b>0$, but one can also consider the case $a=0$, which makes the nonlocal term $M$ become degenerate, and models a string with zero tension when undeformed.

Many other physical phenomena can be modeled through nonlocal equations similar to \pref{eq_main_W_f}  (see examples in \cite{Vill_KirchhAplic,DaiHao09_p(x)Lap}), and interesting  mathematical questions also arise.  
For more  recent literature about such  Kirchhoff type problems  we  cite  the works 
 \cite{AlvCorMa,Ma_kirchhSurv,CorFig_Kirchh,ColaPuc_KirchhDeg,Anello_kirch,ChengWuLiu_kirch_conc,TaCh16_kirchh_GtStSiChg,SoChTa16_kirchreshigh,AmbArc_Kirch_Bif,AmbArc_Kirch_VarDeg,SanSic18_KircDeg,itumas_contrex,MadNun_kirchCC},   which deal with the existence of  solutions with various types of nonlinearities $f$ and use mainly  variational methods. Among them, we refer to \cite{ColaPuc_KirchhDeg,SoChTa16_kirchreshigh,AmbArc_Kirch_VarDeg,SanSic18_KircDeg,MadNun_kirchCC} for considering also the case where the nonlocal term $M$ is degenerate.
  
We remand to Section \tref{sec_commap} for a further discussion of the existing literature concerning the applications considered in Section \tref{sec_appl}.
 
\par\medskip

This paper is structured as follows: in the sections \tref{sec_nomin} and \tref{sec_simin} we state and prove the main theorems. In Section \tref{sec_PS}, we prove a compactness condition for the functional \pref{eq_funct_gen} in a rather general form, that will be then used in the following Section \tref{sec_appl} to obtain existence of solutions for some nonlocal problem as \pref{eq_main_W_f}.
Section \tref{sec_CE_Lieb} is instead devoted to the mentioned  counterexample concerning regularity and a-priori estimates for degenerate Kirchhoff type problems as  \pref{eq_main_W_f}. Finally, in the Appendix we summarize some estimates, which are used in the proof of Theorem \tref{th_nomin_gen}, involving suitable compact support approximations of the instanton functions which realize the best  constant for the  Sobolev embedding $\W\subseteq L^{\frac{pN}{N-p}}$.

\subsection{Notation and first remarks}\label{sec_NotRem}
Throughout the paper we will denote by  $\n{\cdot}_W=\pt{\int_\Omega|\nabla u|^p}^{1/p}$ the norm in $W_0^{1,p}(\Omega)$ and by $\n{\cdot}_s$
the $L^s$-norm.
As usual,  $p^*=\frac{pN}{N-p}$ will denote the critical exponent of the embedding $W_0^{1,p}(\Omega)\hookrightarrow L^s(\Omega)$, when $N>s$.
We will also use the letters $C,c$ to denote  generic positive constants which may
vary from line to line.

In order to have  the functional \pref{eq_funct_gen} well defined in $\W$ we will always assume that $M$ is integrable at $0$ and continuous on $(0,\infty)$ and that, for some positive constant $C_0$, 
\begin{itemize}
\item[(\aslabel{$H_0$}) \label{Hp_nomin_sub}] \qquad $F(x,v)\leq C_0(1+|v|^{p^*})$, \qquad \text{ for every $x\in\Omega$, $v\in\R$}\,.
\end{itemize}

 Observe that a solution of the nonlocal problem \pref{eq_main_W_f} is   also a solutions of the local problem
\begin{equation}\label{eq_main_gam}
\left\{
\begin{array}{ll}
-\Delta_p u=\llambda f(x,u)& \mbox{ in }\ \Omega\,,
\\
u=0 & \mbox{ on }\ \partial \Omega\,,
\end{array}
\right.
\end{equation}
with $\llambda=1/M(\n{u}_W^{p})$: this observation will be used several times in order to translate known properties of the $p$-Laplacian to the nonlocal case.

\section{H\"older minimizers which are not Sobolev minimizers}\label{sec_nomin}
In this section we show that for certain nonlocal degenerate problems no result similar to \cite{H1versusC1, GPM_sob_vs_hold,GuoZhang_sob_vs_hold,BrItUb_multplap} may hold true: in fact, we show that it is possible to have a minimum with respect to the $L^\infty$ norm (which  then is also a minimum with respect to the $\Cont^1$ norm) which is not  a minimum in $W^{1,p}_0$.

Our result is contained in the following theorem.
\begin{theorem}\label{th_nomin_gen}
Let  $\Omega$ be a bounded domain in $\R^N$ with $N>p$ and consider the functional $J$ as in 
\pref{eq_funct_gen}, where $F$  satisfies hypothesis \pref{Hp_nomin_sub}.
Let
\begin{equation}\label{eq_relexp_nomin}r>p^*\geq\p>q\geq1
\end{equation} and assume that for suitable  constants $C_1,C_2,C_3>0$,
\begin{itemize}
\item[(\aslabel{$P_M^<$}) \label{Hp_nomin_Mle}]  \quad $\frac1p\widehat M(s^p)\leq C_1 s^{r}$,\qquad  for  $s$ small, 
\item[(\aslabel{$P_F$}) \label{Hp_nomin_Fge}]  \quad   $F(x,v)\geq C_2v^\p- C_3v^{q}$,\qquad   for $x\in\Omega,\, v\geq0$. 
\end{itemize}
Then  there exists a sequence $u_n$ in $W^{1,p}_0(\Omega)$ with $\n{u_n}_W\to0$, such that   $J(u_n)<0$.

In particular, if also
\begin{itemize}
\item[(\aslabel{$P_M^0$}) \label{Hp_nomin_Mge0}] \quad   $\widehat M(t)\geq0$, \qquad for $t$ small, 
\item[(\aslabel{$P_F^0$}) \label{Hp_nomin_Fle0}] \quad  $F(x,v)\leq0$,\qquad  for $x\in\Omega$ and $|v|$ small, 
\end{itemize}
then  the origin is a local minimum for $J$ with respect to the $\Cont^1$   norm, but  not with respect to the  $W_0^{1,p}$ norm.

Finally, if 
\begin{itemize}
\item[(\aslabel{$P_M^*$}) \label{Hp_nomin_Mge0all}] \quad $\widehat M(t)\geq0$,\qquad  for every $t\geq0$,  
\end{itemize}
then   the origin is a local minimum for $J$ 
with respect to the $\Cont^1$  and $L^\infty$ norms, but  not with respect to the  $W_0^{1,p}$ norm.
\end{theorem}
\begin{remark}
Observe that for hypotheses \pref{Hp_nomin_Mle} and \pref{Hp_nomin_Mge0} to hold  it is sufficient (but not necessary) that   $ 0\leq M(s^p)\leq C s^{r-p}$ for small $s$.\eor
\end{remark}

In our model examples (\ref{eq_funct_p}\,-\,\ref{eq_funct_p+}), it is straightforward to verify that if \pref{eq_relexp_nomin} holds true and $\lambda>0$, then the hypotheses \pref{Hp_nomin_sub}, \pref{Hp_nomin_Mle}, \pref{Hp_nomin_Fge}, \pref{Hp_nomin_Mge0all} and \pref {Hp_nomin_Fle0} are satisfied and then we have the following result.
\begin{corollary}\label{coro_nomin_model}
If \pref{eq_relexp_nomin} holds true and $\lambda>0$,  then the origin is a local minimum for the functionals   (\ref{eq_funct_p}\,-\,\ref{eq_funct_p+}) 
with respect to the $L^\infty$ norm (and to the $\Cont^1$ norm), but  not with respect to the  $\W$ norm.
\end{corollary}
\begin{proof}[Proof of Theorem \tref{th_nomin_gen}]
By \pref{eq_relexp_nomin}, the following inequalities hold: $$0\leq\frac{N-p}p\frac{p^*-\p}{r-\p} <\frac{N-p}p\,,$$
as a consequence, we can  define the constants $\s,\ttau$, whose role will be clarified later, as below: 
  $$\frac{N-p}p>\s>\frac{N-p}p\frac{p^*-\p}{r-\p}\geq0,$$ 
 $$0<\ttau<(\p-q)\pt{\frac{N-p}p-\sigma}.$$ 

Now let $\pg{\psi_\varepsilon,\ \varepsilon>0}$ be the compact support approximations of the instanton  functions as defined in \pref{eq_defpsi_p}, where the parameter $\beta$ is chosen in such a way that, as  $\varepsilon\to0$, we have, from equations \pref{eq_estPsi_H_p} and \pref{eq_estPsi_s_p} in Corollary \tref{cor_GaRuCa_p},
\begin{eqnarray}
\n{\psi_\varepsilon}_W^r&\leq& C(1+o(1))\leq 2C,\\
\n{\psi_\varepsilon}_q^q&\leq &C\varepsilon^{(N-\frac{N-p}{p}q-\ttau)}\,,\\
\label{eq_nomin_estp}\n{\psi_\varepsilon}_\p^\p&\geq& c\varepsilon^{(N-\frac{N-p}{p}\p)},
\end{eqnarray}
(here the constants $C,c>0 $ and $\beta\in(0,1)$  depend on $N,p,q,\p$, which are fixed along this proof). Observe that \pref{eq_nomin_estp} becomes $\n{\psi_\varepsilon}_\p^\p\geq c$  when $\p=p^*$, which is also true by \pref{eq_estPsi_p*_p}.

By hypotheses \pref{Hp_nomin_Mle} and \pref{Hp_nomin_Fge}  (observe that $\psi_\varepsilon\geq0$) we have 
$$\frac1p\widehat M (\n{u}^p_W)\leq C_1 \n{u}^r_W,\qquad \text{for small $\n{u}_W$},$$
 $$\int _\Omega F(x,u)\geq C_2\n{u}_\p^\p- C_3\n{u}_q^{q},$$
then for $t>0$ small enough we get
\begin{eqnarray}
\nonumber J(t\psi_\varepsilon)&\leq&C_1\n{t\psi_\varepsilon}_W^r+C_3\n{t\psi_\varepsilon}_q^q-C_2\n{t\psi_\varepsilon}_\p^\p\\
\nonumber &=&C_1{t^r}\n{\psi_\varepsilon}_W^r+C_3{t^q}\n{\psi_\varepsilon}_q^q-C_2{t^\p}\n{\psi_\varepsilon}_\p^\p\\
\label{eq_Jteps} &\leq& C(t^r+t^q\varepsilon^{N-q\frac{N-p}p-\ttau})-ct^\p\varepsilon^{N-\p\frac{N-p}p}\,.
\end{eqnarray}

Now we let $t$ depend on $\varepsilon$ by the relation $t=\varepsilon^\s $, so that
\begin{equation}\label{eq_nomin}
J(\varepsilon^\s\psi_\varepsilon)\leq C\pt{\varepsilon^{r\s }+\varepsilon^{N+q\pt{\s -\frac{N-p}p}-\ttau}}-c\varepsilon^{N+\p\pt{\s -\frac{N-p}p}}\,.
\end{equation}
Since $\s>\frac{N-p}p\frac{p^*-\p}{r-\p}\geq0$, the 
first exponent in  \pref{eq_nomin} is  larger than the last  one, moreover, since $\s <\frac{N-p}p$ and $\ttau<(\p-q)\pt{\frac{N-p}p-\sigma}$,  the second exponent is also   larger than the last  one. 

We conclude that  $\n{\varepsilon^\s  \psi_\varepsilon}_W\to0$  as $\varepsilon\to 0$, but  for $\varepsilon$ small enough, $J(\varepsilon^\s  \psi_\varepsilon)<0=J(0)$ because the dominant term in \pref{eq_nomin} is the last one.
A sequence as in the claim would be for instance $u_n=(1/n)^\s  \psi_{1/n}$.

Finally, if $u\in \Cont^1(\Omega)$ 
 with $\n{u}_{\Cont^1}$ small, then  also the $\W$ norm is small and by (\ref{Hp_nomin_Mge0}\,-\,\ref{Hp_nomin_Fle0})  we obtain $J(u)\geq 0$, that is, the origin is a local minimum with respect to the $\Cont^1$ norm.
If also \pref{Hp_nomin_Mge0all} holds  then $J(u)\geq 0$ also for  $u\in L^\infty(\Omega)\cap W^{1,p}_0(\Omega)$ 
 with $\n{u}_{L^\infty}$ small, that is, the origin is a local minimum also with respect to the $L^\infty$ norm.
\end{proof}

\begin{remark}
Of course,   $\varepsilon^\s  \psi_\varepsilon$ is not bounded in $L^\infty$ as $\varepsilon\to0$, in fact, by \pref{eq_def_phi} we get  $\varepsilon^\s\psi_\varepsilon(0)\simeq\varepsilon^{\s+(p-N)/p}\to\infty$.
However, when  (\ref{Hp_nomin_Mge0}\,-\,\ref{Hp_nomin_Fle0}) hold true, along the segment from the origin to  $\varepsilon^\s  \psi_\varepsilon$, $J$ is initially positive, before reaching its negative value at $\varepsilon^\s  \psi_\varepsilon$.
\eor
\end{remark}

\section{H\"older minimizers which are also  Sobolev minimizers}\label{sec_simin}
Our purpose in this section is to show that when the degeneration of the function $\widehat M$ in the functional \pref{eq_funct_gen} is weaker than in Theorem \tref{th_nomin_gen}, it is instead possible to prove that   minima with respect to the $L^\infty$ norm are also minima  with respect to the $W^{1,p}_0$ norm. Under  more restrictive conditions, we can also obtain the stronger result that  minima with respect to the $\Cont^1$ norm are also minima  with respect to the $W^{1,p}_0$ norm, as in the  classical results.

The first result is the following.
\begin{theorem}\label{th_simin} 
 Let  $\Omega$ be a  bounded and smooth  domain in $\R^N$ with $N>p$ and consider the functional $J$ as in 
\pref{eq_funct_gen}, associated to the nonlocal equation \pref{eq_main_W_f}.
Suppose the function  $f:\Omega \times \R\to\R$ in \pref{eq_main_W_f} is a continuous function such that  \pref{Hp_nomin_sub} holds true and 
\begin{itemize}
\item[(\aslabel{$Q_f$}) \label{fest1}]\qquad there exist constants $D>0$ and $\l\in[p,p^*)$ such that 
$$\displaystyle f(x,v)sgn(v)\leq D|v|^{\l-1},\qquad\forall\,(x,v)\in\Omega\times\R.$$
\end{itemize}
Assume further that 
\begin{itemize}
\item[(\aslabel{$Q_M^>$}) \label{Hp_min_Morig}] 
\qquad $M(t)\geq0$ for every $t\geq0$ and  there exist constants $a_1,\delta>0$ and $r\in(p,p^*)$ such that
$$M(s^p)\geq\frac{r\,a_1}{p}s^{{r}-p}, \qquad\text{ for $0\leq s^p<\delta$. }$$ 
\end{itemize} 

Then, if the origin is a local minimum for $J$ with respect to  the $L^\infty$ norm, then it is also a local minimum with respect to the $W_0 ^{1,p}$ norm.
\end{theorem}

\begin{remark}
The hypothesis \pref{Hp_min_Morig} in Theorem \tref{th_simin} implies that $\widehat M(s^p)\geq a_1s^r$ for $s$ small, then it is analogous (but with reversed inequality) to hypothesis \pref{Hp_nomin_Mle} in Theorem \tref{th_nomin_gen}. The condition here, however, is given on $M$ instead of its primitive, then it is slightly stronger. \eor
\end{remark}

For the model functionals \pref{eq_funct_p}-\pref{eq_funct_p+} the hypotheses of Theorem 
\tref{th_simin} are satisfied provided 
\begin{equation}\label{eq_relexp_simin}
0<r<p^*
\qquad and \qquad
p^*>\p>q\geq1,
\end{equation}
moreover, the origin is a local minimum with respect to the $L^\infty$ norm since (\ref{Hp_nomin_Mge0all}\,-\,\ref{Hp_nomin_Fle0}) hold true. 
Then we have the following result.
\begin{corollary}\label{coro_simin_model}
 Let  $\Omega$ be a  bounded and smooth  domain in $\R^N$ with $N>p$.    If \pref{eq_relexp_simin} holds true, then  the origin is a local minimum for   (\ref{eq_funct_p}) and   (\ref{eq_funct_p+}) with respect to the  $L^\infty$ norm, and  also with respect to the $W_0^{1,p}$ norm.  
\end{corollary}

In order to prove Theorem \tref{th_simin} we will need the a-priori estimate contained in the following  lemma, which is obtained through the  application of  Moser's iterations.

\begin{lemma}\label{lm_estim} Let  $\Omega$ be a  bounded and smooth  domain in $\R^N$ with $N>p$
and  $f:\Omega \times \R\to\R$ be a continuous function which
satisfies hypothesis  \pref{Hp_nomin_sub} and 
\begin{equation}\label{fest1_lm}
f(x,v)sgn(v)\leq L|v|^{\q-1}\qquad\forall\,(x,v)\in\Omega\times\R,
\end{equation}
for some numbers $L>0$ and $\q\in[p,p^*)$.

Then there exists  $C_1(\q,p,\Omega) $ such that 
\begin{equation}\label{eq_est_Moser}
\n{u}_\infty\leq  C_1(\q,p,\Omega) L^{\frac1{p^*-\q}} \n{u}_{p^*}^{\frac{p^*-p}{p^*-\q}}\,,
\end{equation} 
for every weak solution  $u\in
W_0 ^{1,p} (\Omega )$   of the problem
\begin{equation}\label{eq_inMoser}
\left\{
\begin{array}{ll}
-\Delta _p u= f(x,u) & \mbox{ in }\ \Omega\,,
\\
u=0 & \mbox{ on }\ \partial \Omega\,.
\end{array}
\right. 
\end{equation}

\end{lemma}
\begin{proof}
The proof is the same as in \cite[Lemma 2.2]{itulormon}, one only has to observe that $f$ can depend on $x$ and (\ref{fest1_lm}\,-\,\ref{Hp_nomin_sub}) can substitute the stronger condition $|f(v)|\leq L|v|^{\q-1}$. 
  Moreover, observe that  \pref{Hp_nomin_sub} is used only to guarantee that the solutions are  $C^{1,\alpha}(\overline\Omega)$, so that the final estimate depends on $\q$ and $L$ but not on $C_0$.
\end{proof}
\begin{remark}
 Lemma \tref{lm_estim},  and then Theorem \tref{th_simin}, can be extended to the case $N=p$ by replacing $p^*$ with any large number. 
If instead $N<p$ then Theorem \tref{th_simin} holds true trivially, since $\W\subseteq L^\infty$.
\eor
\end{remark}

\begin{proof}[Proof of Theorem \tref{th_simin}]
Our proof is based on the proof of \cite[Lemma 2.2]{BrItUb_multplap}.
In order to be able to exploit Lemma \tref{lm_estim}, we will compare the problem \pref{eq_main_W_f}  with its local versions \pref{eq_main_gam} with $\llambda=1/M(\n{u}_W^p)$.
Observe that \pref{Hp_min_Morig}  implies that $\widehat M$ is nondecreasing and
\begin{equation}\label{eq_estM}
\begin{cases}
 \widehat M(t)\geq a_1{t^{r/p}}\qquad&\text{ for }0\leq t\leq\delta,\\
\widehat M(t)\geq a_1{\delta^{r/p}}\qquad&\text{ for }t\geq \delta.
\end{cases}
\end{equation}

For sake of contradiction, we assume from now on that   the origin 
is not a local minimum for $J$ in the $W_0^{1,p}$ topology, and we aim to prove that then it cannot be a minimum with respect to the $L^\infty$ norm. 

We start by defining
$$
Q(u) : = \frac{1}{\q}\int_{\Omega } |u(x)|^{\q},
\qquad \text{ for $u\in W_0 ^{1,p} (\Omega )$ }
$$
and, for $\varepsilon >0$,
$$
B^Q_{\varepsilon} := \{ u\in W_0 ^{1,p} (\Omega ):\, Q(u ) \leq
\varepsilon \} \,.
$$
First we claim that 
 for every $\varepsilon>0 $ small enough  there exists $v_{\varepsilon }\in B^Q_{\varepsilon }$ such that 
$J
(v_\varepsilon)=\min_{u\in B^Q_\varepsilon}J(u)<0$.
Actually, since $J(0)=0$ and $B^Q_\varepsilon$ contains a small ball in the $W_0^{1,p}$  norm, we get that  $\inf_{u\in B^Q_\varepsilon}J(u)<0$ by our contradiction assumption. 
Then, given  a minimizing sequence $w_n$ for $J$ in $B^Q_\varepsilon$, one may assume 
$\frac1p\widehat M(\n{w_n}_W^p)\leq \int F(x,w_n)$; 
 by \pref{fest1}, this can be controlled by $\varepsilon$ so, for  $\varepsilon>0 $ small enough, we may assume  $\widehat M(\n{w_n}_W^p)< a_1{\delta^{r/p}}$  and then, using \pref{eq_estM}, $\n{w_n}_W\leq\delta^{1/p}$.
 As a consequence,  we may assume that $w_n$ converges weakly in $\W$ and strongly in $L^\q$ to some $v_\varepsilon\in B^Q_\varepsilon$ with $J(v_\varepsilon)= \min_{B^Q_\varepsilon} J<0$ (observe that since $\widehat M$ is nondecreasing $J$ is weakly lower semicontinuous).

Now since, as $\varepsilon\to0$, $v_\varepsilon\to 0$ in $L^\q$ and  $J(v_\varepsilon)<0$, reasoning as above we deduce that  $\widehat M(\n{v_\varepsilon}_W^p)\to0$ and then  $v_\varepsilon\to0$ in $W_0^{1,p}$.
We consider two cases. 
\medskip

\noindent {\bf Case 1}:  $Q(v_{\varepsilon})<\varepsilon $. Then
$v_{\varepsilon }$ is also  a local minimizer of $J$ in $
W^{1,p}_0(\Omega)$ (actually, there exists a small ball centered at $v_\varepsilon$ and contained in $B^Q_\varepsilon$).
 Hence $v_{\varepsilon} $ is a solution of
\pref{eq_main_W_f}, and also of \pref{eq_main_gam} with $\llambda_\varepsilon= M(\n{v_\varepsilon}_W^p)^{-1}$.
By hypotheses (\ref{fest1}\,-\,\ref{Hp_nomin_sub})  we may apply Lemma \tref{lm_estim} to equation \pref{eq_main_gam} with $L:=D \llambda_\varepsilon$ to obtain 
 $$\n{v_\varepsilon}_\infty\leq  C_1(\q,p,\Omega) \pt{D\llambda_\varepsilon}^{1/(p^*-\q)} \n{v_\varepsilon}_{p^*}^{(p^*-p)/(p^*-\q)}\,,$$  
 where since $\llambda_\varepsilon= M(\n{v_\varepsilon}_W^p)^{-1}\leq \frac{p}{ra_1}\n{v_\varepsilon}_W^{p-r}$ we get 
 \begin{eqnarray}\label{eq_est_vep_W}
\nonumber\n{v_\varepsilon}_\infty&\leq& C_1(\q,p,\Omega)\pt{\frac{p}{ra_1}D}^{1/(p^*-\q)}\n{v_\varepsilon}_W^{(p-r)/(p^*-\q)}\n{v_\varepsilon}_{p^*}^{(p^*-p)/(p^*-\q)}
 \\\label{eq_stimLinfty}&\leq&  C(\q,p,\Omega)\n{v_\varepsilon}_W^{\frac{p^*-r}{p^*-\q}};
\end{eqnarray}
 since  $r<p^*$, we get that the $L^\infty$ norm goes to zero as $\varepsilon\to0$. 

We have thus proved that the origin is not a local minimum with respect to the $L^\infty$ norm.

\medskip

\noindent {\bf Case 2}:  $ Q(v_{\varepsilon} ) =\varepsilon $. Then
there exists a Lagrange
multiplier $\mu _{\varepsilon } \in \R $ such that $J'(v_\varepsilon) = \mu_\varepsilon
Q'(v_\varepsilon)$, i.e., $v_{\varepsilon} $ is a weak solution of
\begin{equation}
\label{eqveps} \left\{
\begin{array}{ll}
-M(\n{v_{\varepsilon}}_W^p)\Delta _p  v_{\varepsilon} = f(x,v_{\varepsilon})+\mu_{\varepsilon} |v_{\varepsilon} | ^{\q-2} v_{\varepsilon}
& \mbox{ in } \  \Omega ,
\\
u=0 & \mbox{ on } \ \partial \Omega ,
\end{array}
\right.
\end{equation}
and a weak solution of
\begin{equation}
\label{eqveps_lm} \left\{
\begin{array}{ll}
-\Delta _p v_{\varepsilon} =\llambda_\varepsilon  \pq{f(x,v_{\varepsilon})+\mu_{\varepsilon} |v_{\varepsilon} | ^{\q-2} v_{\varepsilon}}
& \mbox{ in } \  \Omega ,
\\
u=0 & \mbox{ on } \ \partial \Omega ,
\end{array}
\right.
\end{equation}
  with $\llambda_\varepsilon=M(\n{v_\varepsilon}_W^p)^{-1}$.
By testing equation \pref{eqveps} with $v_{\varepsilon }$ one gets 
 $${J'(v_\varepsilon)[v_\varepsilon]} = \mu_\varepsilon {Q'(v_\varepsilon)[v_\varepsilon]}=\mu_\varepsilon\int |v_\varepsilon|^\q= \mu_\varepsilon \q\,\varepsilon\,.$$
If $\mu_\varepsilon>0$ then ${J'(v_\varepsilon)[v_\varepsilon]}>0$ which implies that $J$  decreases in the direction $-v_\varepsilon$, contradicting that $v_\varepsilon$ is a global minimizer for $J$ in $B^Q_\varepsilon$.
It follows that $\mu_{\varepsilon} \leq 0 $. As a consequence, the right hand side  of \pref{eqveps_lm} satisfies the hypotheses of Lemma \tref{lm_estim} with $L=D \llambda_\varepsilon$, which allows us to  obtain again \pref{eq_stimLinfty} and to get  as before that the origin is not a minimum with respect to the $L^\infty$ norm.
\end{proof}

The second result in this section aims to provide a version of the more classical result that relates minima with respect to the $\Cont^1$ norm with those in the Sobolev norm. 
\begin{theorem}\label{th_simin_C1} 
In the hypotheses of Theorem \tref{th_simin}, if moreover 
\begin{equation}\label{eq_condr_C1H1}
(\l-1)> \, (r-p)\, \frac{p^*-1}{p^*-p}
\end{equation}
and there exist further constants  $\widetilde D,\lt>0$  such that 
 \begin{equation}\label{eq_condlt_C1H1}
\begin{cases}\displaystyle |f(x,v)|\leq \widetilde D|v|^{\lt-1},\qquad\text{for $x\in\Omega$, $|v|$ small,}\\ \displaystyle\lt-1>(r-p)\frac{(p^*-\q)}{(p^*-r)}\,,
\end{cases}
\end{equation}
then, if the origin is a local minimum for $J$ with respect to  the $\Cont^1$  norm, then it is also a local minimum with respect to the $W_0 ^{1,p}$ norm.
\end{theorem}
\begin{remark}\label{rmk_l_lt}
The conditions (\ref{eq_condr_C1H1}\,-\,\ref{eq_condlt_C1H1}) state a balance between the growth of the nonlocal term near zero and the behavior of the nonlinearity $f$. In fact, a higher degeneration of the nonlocal term needs to be compensated with a higher exponent in \pref{fest1} and in \pref{eq_condlt_C1H1}, if one aims to control the $\Cont^1$ norm.
\\
In the case of the model functionals \pref{eq_funct_p}-\pref{eq_funct_p+} satisfying \pref{eq_relexp_simin}, the exponents $\l$ and $\lt$ correspond, respectively, to $\p$ and $q$. 
Observe that, since $r<p^*$, there always exists  $\l<p^*$ satisfying  \pref{eq_condr_C1H1}. If $\l$ is close to $p^*$, then condition \pref{eq_condlt_C1H1} is not very restrictive, allowing the nonlinearity to have a lousier estimate from below near zero.
\eor
\end{remark}

\begin{proof}[Proof of Theorem \tref{th_simin_C1}]
As in the proof of Theorem \tref{th_simin}, we assume for sake of contradiction  that   the origin 
is not a local minimum for $J$ in the $W_0^{1,p}$ topology, and we aim to prove now that it cannot  be a minimum with respect to the $\Cont^1$ norm.

Carrying on from the proof of Theorem \tref{th_simin}, the first step will be  to prove that if  (\ref{eq_condr_C1H1}\,-\,\ref{eq_condlt_C1H1}) hold true, then it is possible to  obtain a uniform  bound, as $\varepsilon\to0$, for $\n{v_\varepsilon}_{\Cont^{1,\alpha}}$, for some $\alpha>0$.
We already proved  that \pref{eq_stimLinfty} holds true and then $\n{v_\varepsilon}_W,\,\n{v_\varepsilon}_{\infty}\to0$.

In case 1), since  $v_\varepsilon$ satisfies \pref{eq_main_gam} with $\llambda_\varepsilon= M(\n{v_\varepsilon}_W^p)^{-1}$, it also   satisfies
\begin{equation}\label{fproblem_tr}
\left\{
\begin{array}{ll}
-\Delta _p u= \llambda_\varepsilon f_\varepsilon(x,u) & \mbox{ in }\ \Omega\,,
\\
u=0 & \mbox{ on }\ \partial \Omega\,,
\end{array}
\right. 
\end{equation}
where $f_\varepsilon(x,u)$ is truncated such that  $f_\varepsilon(x,t)=f(x,sgn(t)\n{v_\varepsilon}_\infty)$ for $|t|>\n{v_\varepsilon}_\infty$, implying that, for $\varepsilon$ small, $|f_\varepsilon(x,t)|\leq \widetilde D\n{v_\varepsilon}_\infty^{\lt-1}$, by \pref{eq_condlt_C1H1}.
Then the right hand side of \pref{fproblem_tr} can be estimated as 
\begin{equation}\label{eq_rhs<n}
 |\llambda_\varepsilon f_\varepsilon(x,u)|\leq C \n{v_\varepsilon}_W^{\frac{(p^*-r)}{(p^*-\q)}(\lt-1)-r+p}.
\end{equation}
The exponent in \pref{eq_rhs<n} is positive by  \pref{eq_condlt_C1H1}, providing a uniform bound as  $\varepsilon\to0$. We can thus apply the results in \cite[Theorem 1]{Lieberman} to obtain the uniform  bound for $\n{v_\varepsilon}_{\Cont^{1,\alpha}}$ and some $\alpha>0 $. 

In case 2), we already proved that $v_\varepsilon$ satisfies \pref{eqveps_lm}  with $\mu_\varepsilon\leq0$; 
by testing this equation with $v_\varepsilon$ we have 
$$0\leq \int |\nabla v_\varepsilon|^p=\gamma_\varepsilon \int \pq{f(x,v_\varepsilon)v_\varepsilon+\mu_\varepsilon v_\varepsilon^\q}\,,$$
which implies, using \pref{fest1}, that
 $$-\mu_\varepsilon \int v_\varepsilon^\q\leq\int f(x,v_\varepsilon)v_\varepsilon\leq  D\int v_\varepsilon^\q\,,$$
so that $\mu_\varepsilon\in[-D,0]$.
With this information we may proceed as in case 1), truncating both $f$ and the additional term $\mu_{\varepsilon} |v_{\varepsilon} | ^{\q-2} v_{\varepsilon}$ in \pref{eqveps_lm},  whose right hand side is then estimated as
\begin{equation}\label{eq_rhs<n_mu}
 \m{\llambda_\varepsilon\pq{f(x,v_{\varepsilon})+\mu_{\varepsilon} |v_{\varepsilon} | ^{\q-2} v_{\varepsilon}}}\leq C \pt{\n{v_\varepsilon}_W^{\frac{(p^*-r)}{(p^*-\q)}(\lt-1)-r+p}+\n{v_\varepsilon}_W^{\frac{(p^*-r)}{(p^*-\q)}(\l-1)-r+p}} .
\end{equation}
Using again  \pref{eq_condlt_C1H1}  and now also \pref{eq_condr_C1H1}, which is equivalent to 
${\frac{(p^*-r)}{(p^*-\q)}(\l-1)-r+p}>0$, we obtain again the uniform bound on  $\n{v_\varepsilon}_{\Cont^{1,\alpha}}$, via \cite{Lieberman}.

Finally, with this uniform bound, we can  apply  Ascoli-Arzela's theorem to obtain  a sequence $\varepsilon_n\to 0$  such that
$v_{\varepsilon_n} \to 0$ in $\Cont^1$. 
We have thus proved that the origin cannot be a minimum with respect to the $\Cont^1$ norm.
\end{proof}

\section{Compactness condition}\label{sec_PS}
In this section we will provide conditions under which the functional \pref{eq_funct_gen}  satisfies the compactness conditions required in order to apply classical variational methods. 
\begin{proposition}\label{lm_PSgen}
Let $\Omega$ be a bounded domain and consider the functional \pref{eq_funct_gen}, associated to  problem \pref{eq_main_W_f}. 
  Suppose 
\begin{itemize}
\item[(\aslabel{$\widetilde H_1$})\label{Hp_PS_cresc}] \qquad 
 $f$ is a Carathéodory function   and there exists $\widetilde\p\in(1,p^*)$ and $\widetilde C_0>0$ such that
 $$|f(x,v)|\leq \widetilde C_0(1+|v|^{\widetilde\p-1})\qquad \text{ for every $x\in\Omega$, $v\in\R$}\,,$$
\end{itemize}
 and one of the following two conditions:
\begin{itemize}
\item[(\aslabel{$K_{AR}$}) \label{Hp_PS_SQ}] \qquad 
  there exist $\widetilde r ,C,\beta>0$ such that  
\begin{eqnarray*}
&&(i)\ \qquad \nowlabel{\ref{Hp_PS_SQ}-i} \label{Hp_PS_SQi} 
\widetilde r  F(x,v)-f(x,v)v\leq C\qquad \text{ for every  $x\in\Omega$, $v\in\R$},\\
&&(ii)\qquad \nowlabel{\ref{Hp_PS_SQ}-ii} \label{Hp_PS_SQii}
 \frac{\widetilde r }p\widehat M(s^p)-M(s^p)s^p\geq \beta s-C\qquad \text{ for every $s\geq0$,}
\end{eqnarray*}
\item[(\aslabel{$K_{C}$}) \label{Hp_PS_sQ}]\qquad 
 there exist   $\widetilde r >\widetilde\p$, $C>0$ such that 
$$\widehat M(s^p)\geq s^{\widetilde r }-C\qquad \text{ for every $s\geq0$}.$$
\end{itemize}
 Then 
any PS sequence for the functional \pref{eq_funct_gen}  is bounded.

If \pref{Hp_PS_SQii} is replaced by the weaker condition 
\begin{itemize}
\item[(\aslabel{\ref{Hp_PS_SQ}-iii}) \label{Hp_PS_SQiii}]\qquad there exists $h>0$ such that 
$$ \qquad \frac{\widetilde r }p\widehat M(s^p)-M(s^p)s^p\geq \beta s^h-C\qquad \text{ for every $s\geq0$},$$  
\end{itemize}
 then 
any Cerami sequences for the functional \pref{eq_funct_gen}  is bounded.

The PS (resp Cerami) condition is satisfied provided we assume also
\begin{itemize}
\item[(\aslabel{$K_M$}) \label{Hp_PS_Mpos}]  \qquad $M(t)>0$ and is continuous, for $t>0$.  
\end{itemize}
\end{proposition}
\begin{remark}
The conditions assumed in Proposition \tref{lm_PSgen} are similar to those used in \cite{MadNun_kirchCC} and, for p=2, also in  \cite{AmbArc_Kirch_Bif,AmbArc_Kirch_VarDeg}.

Condition \pref{Hp_PS_SQ} is a generalization of the classical Ambrosetti-Rabinowitz superlinearity condition. Actually, in the local case $M\equiv1$, 
condition \pref{Hp_PS_SQii} reads
$$\pt{ \frac{\widetilde r }p-1} s^p\geq \beta s-C$$ and then reduces to the usual requirement that $\widetilde r >p$ in  \pref{Hp_PS_SQi}. More in general, \pref{Hp_PS_SQ} establishes that ``$f$ grows more than $M$'' at infinity, while \pref{Hp_PS_sQ} establishes the opposite.\eor
\end{remark}
\begin{proof}[Proof of Proposition \tref{lm_PSgen}]
If $u_n$ is a PS sequence for $J$ then one has, for some sequence $\varepsilon_n\to0$,
$$\m{J(u_n)}=\m{\frac1p  \widehat M\left( \n{ u_n}^p_W\right)-\int_{\Omega }
F(x,u_n)}\leq C,$$
$$\m{J'(u_n)[u_n]}=\m{ M\left( \n{ u_n}^p_W\right) \n{ u_n}^p_W-\int_{\Omega }
f(x,u_n)u_n}\leq \varepsilon_n\n {u_n}_W.$$
If condition \pref{Hp_PS_SQ} holds, then one can estimate $\m{\widetilde r  J(u_n)-J'(u_n)[u_n]}$, obtaining
\begin{multline*}\m{\pt{\frac{\widetilde r }p \widehat M\left( \n{ u_n}^p_W\right) - M\left( \n{ u_n}^p_W\right) \n{ u_n}^p_W}-\pt{\int_{\Omega }
\widetilde r  F(x,u_n)-f(x,u_n)u_n}}\leq C+\varepsilon_n\n {u_n}_W,
\end{multline*}
so that 
\begin{equation}\label{eq_unbdd}
\beta \n {u_n}_W \leq C'+\varepsilon_n\n {u_n}_W,
\end{equation} which implies that $\n {u_n}_W$ is bounded. 

If instead condition  \pref{Hp_PS_sQ} is assumed then  $\n {u_n}_W$ is bounded since, in view of   \pref{Hp_PS_cresc}, the functional becomes coercive. 

Finally, if  condition \pref{Hp_PS_SQiii} substitutes \pref{Hp_PS_SQii}  and  $u_n$ is a Cerami sequence, that is, if we assume the stronger condition $(1+\n{u_n})J'(u_n)\to0$, 
then we obtain  $$\beta \n {u_n}_W^h \leq C'+\varepsilon_n$$ instead of \pref{eq_unbdd}, and then again $\n {u_n}_W$ is bounded.

Now by standard arguments it follows that, up to a subsequence, $u_n$ converges weakly in $\W$ to some $u$ and  from  $|J'(u_n)[u_n-u]|\to0$ and \pref{Hp_PS_cresc}, one obtains  $$ M(\n{u_n}^p_W)\int |\nabla u_n|^{p-2}\nabla u_n\nabla (u_n-u)\to0\,.$$
By \pref{Hp_PS_Mpos}
we may assume that either $u_n\to0$ or  $M\pt{\n{u_n}^p_W}\to c>0$. 
In the latter case we have  $\int |\nabla u_n|^{p-2}\nabla u_n\nabla (u_n-u)\to0$ from which, by the $S^+$ property of the $p$-Laplacian, we obtain that $u_n\to u$ strongly in $\W$.
\end{proof}

For   the model functionals (\ref{eq_funct_p}\,-\,\ref{eq_funct_p+}) we have the following.
\begin{corollary}\label{coro_PS_mod}
Let $1<q<\p<p^*$.
\\If  $1\leq r\neq \p$, then 
the functionals (\ref{eq_funct_p}\,-\,\ref{eq_funct_p+})  satisfy the PS condition.
\\If $r\in(0,1)$ then they satisfy the Cerami Condition.
\end{corollary}
\begin{proof}
We only have to check the hypotheses of Proposition  \tref{lm_PSgen}.  
 Condition \pref{Hp_PS_Mpos} always hold true.
Condition \pref{Hp_PS_cresc} is a consequence of the condition $1<q<\p<p^*$.

Condition \pref{Hp_PS_SQ}  is satisfied when $1\leq r<\p$, by taking  $\widetilde r \in(r,\p)$. When $r<1$,
condition \pref{Hp_PS_SQiii} holds  with $h=r$, actually (see \pref{eq_MF_mod}),
$$\frac{\widetilde r }p\widehat M(s^p)-M(s^p)s^p=\pt{\frac{\widetilde r }r-1}s^r\,.$$ 

On the other hand, conditions \pref{Hp_PS_sQ} is satisfied when $r>\p$, by taking $\widetilde\p=\p$ and $\widetilde r \in(\p,r)$.
\end{proof}

\begin{remark}\label{rmk_subterm}
Observe that the term $\frac1q\n{u}_q^q$ has no influence on Corollary \tref{coro_PS_mod}: it can be removed or substituted with a more general term whose growth at infinity is smaller than the  term $\frac1\p\n{u}_\p^\p$. 
\eor
\end{remark}

\section{Applications}\label{sec_appl}
As an application of our main results, we will consider the problem of finding nonnegative  solutions for problem \pref{eq_main_W_f}, in several different situations. 

Nonnegative and nontrivial solutions $u$ can be found as critical points of the modified functional 
\begin{equation} \label{eq_funct_gen+}
J^+(u)=\frac1p\widehat M (\n{u}^p_W)-\int_\Omega F(x,u^+)\,,\qquad u\in W_0^{1,p}(\Omega)\,,
\end{equation}
as can be easily seen by using $u^-$ as a test function. In some cases it will be possible to guarantee that such nonnegative and nontrivial solutions are in fact strictly positive, using the maximum principle in \cite[Theorem 5]{Vazquez} (see Section \tref{sec_commap}).

\subsection{Coercive problems} 
We first consider the model problem \pref{eq_rpq_p}, when the functional associated is coercive, that is, when $r> \p$. We have the following result.
\begin{theorem} \label{th_coe_mod}
Let $1<q<\p<p^*$, $r>\p$ and $\lambda>0$.
\begin{enumerate}
\item If
 $r\in(\p,p^*)$, 
then there exist $\Lambda_2\geq\Lambda_1>0$ such that problem \pref{eq_rpq_p} has no solution (even sign changing) for $\lambda\in(0,\Lambda_1)$ and at least two nonnegative nontrivial solutions for $\lambda>\Lambda_2$.

\item 
If  $r=p^*$, then 
there exist $\widehat \Lambda_2\geq\widehat \Lambda_1>0$ such that problem \pref{eq_rpq_p} has no solution (even sign changing) for $\lambda\in(0,\widehat \Lambda_1)$ and at least one nonnegative nontrivial solutions for $\lambda>\widehat\Lambda_2$.

\item If $r>p^*$, then there exists at least one nonnegative nontrivial solution for every $\lambda>0$. 
\end{enumerate}
\end{theorem}
\begin{proof}[Proof of Theorem \tref{th_coe_mod}]
In the given hypotheses  the functional $\J^+$ in \pref{eq_funct_p+} satisfies the PS condition by Corollary \tref{coro_PS_mod} and is coercive, then there exists a global minimizer.

When  $r>p^*$ we are in the conditions of Theorem \tref{th_nomin_gen} (see Corollary  \tref{coro_nomin_model}) and then  $\J^+(u)<0$ for some $u\in\W$. Then the  global minimum  lies at a negative level and then it corresponds to a nontrivial nonnegative solution.  

If $r\in(\p,p^*)$, on the other hand, we know from Theorem \tref{th_simin}  (see Corollary  \tref{coro_simin_model}) that  the origin is  a local minimum.
In order to prove the existence of two solutions for $\lambda$ large we proceed in the following way:
for  a fixed positive $\phi\in\W$, one has $\J^+(\phi)\to -\infty$ as $\lambda\to \infty$, as a consequence,
 there exists $\Lambda_2>0$ such that for $\lambda>\Lambda_2$,  $i:=\inf_{u\in \W} \J^+(u)<0$. 
We can therefore obtain a solution as a global minimizer for $\J^+$ and a second solution by applying the Mountain Pass theorem (considering paths that connect the origin with the global minimizer).

In the case $r=p^*$ we obtain, by the same argument, that  $\inf_{u\in \W} \J^+(u)<0$ for $\lambda$ large enough, and then we are able to obtain a solution as a global minimizer for $\J^+$.

\par\medskip

The nonexistence result is more difficult to prove. We will follow some ideas from \cite{Ane12_2nl_blw}.
First observe that if $u\in \W$ is a solution of \pref{eq_rpq_p}  then  the functional $\J$ in \pref{eq_funct_p} satisfies $$\J'(u)[u]=\n{u}_W^r+\n{u}_q^q-\lambda\n{u}_\p^\p=0\,.$$
Moreover, by using  H\"older inequality and then  Sobolev inequality, for  $r\leq p^*$,
 $$\n{u}_\p^\p= \int |u|^{q\frac{r-\p}{r-q}}|u|^{r\frac{\p-q}{r-q}}\leq \n{u}_q^{q\frac{r-\p}{r-q}}\n{u}_r^{r\frac{\p-q}{r-q}}\leq  \n{u}_q^{q\frac{r-\p}{r-q}}(C\n{u}_W)^{r\frac{\p-q}{r-q}}\,,$$
from which, by  Young inequality, we obtain 
$$\n{u}_\p^\p\leq \frac{r-\p}{r-q}\n{u}_q^q+\frac{\p-q}{r-q}(C\n{u}_W)^r\,.$$
As a consequence,
$$0=\n{u}_W^r+\n{u}_q^q-\lambda\n{u}_\p^\p\geq\pt{1-\lambda \frac{\p-q}{r-q}C^r}\n{u}_W^r+\pt{1-\lambda \frac{r-\p}{r-q}}\n{u}_q^q.$$
\\ This shows that  $u\equiv 0$ is the unique solution, when  $\lambda$ is so small that both terms in parentheses are positive.
\end{proof}

A more general version of the third point in Theorem \tref{th_coe_mod} is contained in the following theorem.
\begin{theorem} \label{th_coe_gen}
Problem \pref{eq_main_W_f} possesses at least one  nonnegative nontrivial solution  provided  
\begin{itemize}
\item \pref{Hp_PS_cresc}, \pref{Hp_PS_sQ} and   \pref{Hp_PS_Mpos} hold true,
\item  \pref{Hp_nomin_Mle} and \pref{Hp_nomin_Fge} hold true with $r>p^*$ and  $1<q<\p<2^*$.
\end{itemize}
\end{theorem} 
\begin{proof}
The functional $J^+$ in \pref{eq_funct_gen+} satisfies the PS condition by (\ref{Hp_PS_cresc}\,-\,\ref{Hp_PS_sQ}\,-\,\ref{Hp_PS_Mpos}): see Proposition \tref{lm_PSgen}.

Moreover, it is coercive by (\ref{Hp_PS_cresc}\,-\,\ref{Hp_PS_sQ}), however the origin is not a local minimum in view of Theorem \tref{th_nomin_gen}, where we use  \pref{Hp_nomin_Mle} and \pref{Hp_nomin_Fge}  with $r>p^*$ and  $1<q<\p<2^*$.

Then as in the proof of Theorem \tref{th_coe_mod} we have a nonnegative nontrivial solution corresponding to a global minimum at a negative level. 
\end{proof}

\subsection{Non coercive problems}
If $r<\p$ in problem \pref{eq_rpq_p}, then the functional is not coercive anymore, and it is easily seen that a solution exists for every $\lambda>0$,  actually, even dropping the term $-|u|^{q-2}u$, the functional has a minimum at the origin and then a mountain pass structure. 
It turns out to be more interesting then to consider an additional term of intermediate order between $q$ and $\p$: we consider the problem  
\begin{equation}\label{eq_rpqs}
\left\{
\begin{array}{ll}
-\n{u}_W^{r-p}\Delta_p u= -|u|^{q-2}u+\mu |u|^{\s-2}u +\lambda|u|^{\p-2}u & \mbox{ in }\ \Omega\,,
\\
u=0 & \mbox{ on }\ \partial \Omega\,.
\end{array}
\right.
\end{equation}

In our first result  we show that adding this new term,  existence is maintained.
\begin{theorem}\label{th_notcoe_mod}
Suppose   $1<q<r<\p<p^*$ and $\s\in(q,\p)$. Then problem \pref{eq_rpqs} has at least one nonnegative nontrivial solution for all $\lambda>0$ and $\mu\in\R$.
\end{theorem} 
\newcommand{\I}{\mathcal{I}}
\begin{proof}
 Nonnegative solutions will be critical points of the functional 
\begin{equation}\label{eq_funct_rpqs+}
\I^+(u)=\frac1r \n{u}^r_W+\frac1 q\n{u^+}_q^q-\frac\mu\s \n{u^+}_\s^\s-\frac\lambda\p \n{u^+}_\p^\p\,.
\end{equation}

Proceeding as in the proof of Corollary \tref{coro_PS_mod} (see also Remark \tref{rmk_subterm}), one can see that conditions \pref{Hp_PS_cresc},  \pref{Hp_PS_Mpos} and  \pref{Hp_PS_SQ} (with $\widetilde r \in(r,\p)$) are satisfied and then $\I^+$ satisfies the PS condition. 

Moreover, the origin is a local minimum for $\I^+$,  since $r<p^*$, $q<\s<\p$ and we can apply Theorem \tref{th_simin} as in the case of the functional \pref{eq_funct_p+}.

Finally, since $\p$ is the largest power in $\I^+$, there exists $e\in\W$ such that $\I^+(e)<0$.

We can thus apply the Mountain Pass theorem to obtain a nontrivial critical point of $\I^+$.

\end{proof}

If the additional term's exponent $\sigma$ lies below $r$, then we have also the following multiplicity result.
\begin{theorem}\label{th_notcoe_mod_mult}
In the conditions of Theorem \tref{th_notcoe_mod}, if moreover $\sigma\in(q,r)$, then for suitably large $\mu>0$ there exists $\Lambda_\mu>0$ such that for $\lambda\in(0,\Lambda_\mu)$  there exist three nonnegative nontrivial solutions  of problem \pref{eq_rpqs}.
\end{theorem} 

\begin{proof}
From the proof of Theorem \tref{th_notcoe_mod} we know that  the origin is always a local minimum and the PS condition holds true for $\I^+$.
 
Now consider $K(u):=\frac1r \n{u}^r_W+\frac1q \n{u^+}_q^q -\frac \mu\s\n{u^+}_\s^\s.$
 Let   $\phi\in\W$ be a  positive function with $\n\phi_W=1$ and take $S,\mu>0$ (with $\mu$ large enough) such that 
\begin{equation} \label{eq_3sol_0}
\text{$K(\phi)<0\ \ \ $  and $\ \ \  K(\tau \phi)<S\ $ for $\tau\in[0,1]$}.
\end{equation}
 This is possible since $K(\phi)\to-\infty$ as $\mu \to \infty$ and $K$ is bounded in bounded sets.

Having fixed  $\mu$ as above, 
we can find $R>\n\phi_W$ and $\Lambda_\mu>0$ (both depending on $\mu$) such that, for every $\lambda\in(0,\Lambda_\mu)$,
\begin{equation} \label{eq_3sol_1}
\text{$\I^+(u)\geq 2S$ for $\n u_W=R$,}
\end{equation} 
\begin{equation} \label{eq_3sol_2}
\text{ $\I^+(\phi)<0\ \ \ $ and $\ \ \ \I^+(\tau \phi)<S\ $ for $\tau\in[0,1]$}.
\end{equation}
Actually \pref{eq_3sol_1} is possible because 
$$\I^+(u)\geq \frac1r \n{u}^r_W-C\frac \mu\s\n{u^+}_W^\s-C\frac\lambda\p \n{u}_W^\p\,$$ 
and  the sum of the first two terms goes to infinity as $\n{u}_W\to \infty$ (since $\s<r$), while \pref{eq_3sol_2} is a consequence of \pref{eq_3sol_0}.

Finally, since $\p$ is the largest power in the functional, we can also find $T>R$ (depending on $\lambda$) such that  $\I^+(T\phi)<0$. 

With this geometry it is possible to obtain three distinct nontrivial critical points for $\I^+$: 
a local minimum at a negative level inside the ball $\n{u}_W\leq R$, a mountain pass critical point obtained by considering paths that join the local minimum at the origin with $\phi$, whose level is positive but not more than $S$ by \pref{eq_3sol_2}, and a second mountain pass point obtained by considering paths that join the origin with $T\phi$, whose level is at least $2S$ by \pref{eq_3sol_2}.
\end{proof}

 As in the previous section we give here a more general version of the existence result of Theorem \tref{th_notcoe_mod}.
\begin{theorem}\label{th_notcoe_gen}
Problem \pref{eq_main_W_f} possesses a nonnegative nontrivial solution  provided $f$ is a continuous function and 
\begin{itemize}
\item  \pref{Hp_PS_cresc},    \pref{Hp_PS_SQ} and \pref{Hp_PS_Mpos} hold true, 
\item  \pref{Hp_min_Morig} holds with $r<p^*$,
\item[(\aslabel{$H_f$}) \label{Hp_gen1}]  \qquad $f(x,v)sgn(v)\leq0$ \qquad for small $v$,
\item[(\aslabel{$H_B$}) \label{Hp_gen2}] \qquad there exists a ball $B\subseteq\Omega$ such that $$\lim_{v\to+\infty}\frac {F(x,v)}{v^{\widetilde r}}=+\infty\qquad \text{uniformly in $B$},$$
where $\widetilde r$ is the same as in condition  \pref{Hp_PS_SQ}.
\end{itemize}
\end{theorem} 
\begin{proof}
By the conditions \pref{Hp_PS_cresc},  \pref{Hp_PS_Mpos} and  \pref{Hp_PS_SQ}   the functional $J^+$ in \pref{eq_funct_gen+} satisfies the PS condition. 

Since  $M\geq0$ and using  \pref{Hp_gen1}, the origin is a local minimum for $J^+$ with respect to the $L^\infty$ norm. Moreover, conditions \pref{Hp_gen1} and \pref{Hp_PS_cresc} imply that  \pref{fest1} holds true and then  by \pref{Hp_min_Morig}
we can also apply Theorem \tref{th_simin} to obtain that  the origin is  a local minimum. 

Finally, consider a nontrivial function  $\psi\geq 0$, with support in $B$ and $\n\psi_W=1$. 
We will use some estimates that are proved in Lemma \tref{lm_estMF} below:
 let $D$ be the constant in estimate \pref{eq_est_ap2} and take $A=\frac{2D}{p\n\psi^{\widetilde r}_{\widetilde r}}$. By estimate   \pref{eq_est_F} we obtain that 
 $$\int_\Omega F(x,t\psi)\geq A\n{t\psi}^{\widetilde r}_{\widetilde r}-E_A|B|\,,$$
then, using \pref{eq_est_ap2},
$$J^+(t\psi)\leq \frac1p\widehat M(t^p)-  A\n{t\psi}^{\widetilde r}_{\widetilde r}+E_A|B|\leq \frac {D-2D}pt^{\widetilde r}+E_A|B|\,;$$
this proves that $J(t\psi)<0$ for $t$ large enough.

As in the proof of Theorem \tref{th_notcoe_mod} we can thus apply the Mountain Pass theorem to obtain a nontrivial critical point of $J^+$.
\end{proof}
The estimates used in the proof above are contained in the following Lemma.
\begin{lemma}\label{lm_estMF}
Conditions \pref{Hp_PS_Mpos} and \pref{Hp_PS_SQii} imply that there exists $D>0$ such that
\begin{eqnarray}\label{eq_est_ap2}
&&\widehat M(s^p)\leq D s^{\widetilde r},\qquad \text{ for $s$ large enough,} 
\\\label{eq_est_ap3}&&
\lim_{t\to\infty}\widehat M(t)=\infty\,.
\end{eqnarray}
Conditions \pref{Hp_PS_cresc} and \pref{Hp_gen2} imply  that, given $A>0$, there exists $E_A>0$ such that 
\begin{equation}\label{eq_est_F}
F(x,v^+)\geq Av^{\widetilde r}-E_A\qquad  \text {for every $x\in B$.} 
\end{equation} 
\end{lemma}
\begin{proof}
Condition \pref{Hp_PS_SQii} implies that  $\frac{\widetilde r }p\widehat M(s^p)-M(s^p)s^p\geq0$ for $s\geq s_0>0$.
Then rearranging (since everything is positive by \pref{Hp_PS_Mpos}) we obtain  $\frac {M(t)}{\widehat M(t)}\leq \frac{\widetilde r}p \frac1t$ and integrating from $t_0=s_0^{1/p}$ to $t>t_0 $ we get 
$\displaystyle\widehat M(t)\leq\frac{\widehat M(t_0)}{t_0^{\widetilde r/p} }t^{ \widetilde r/p} $, which gives \pref{eq_est_ap2}.
We also have  $\frac{\widetilde r }p\widehat M(s^p)\geq\delta s- C$, which gives \pref{eq_est_ap3}.

By  \pref{Hp_gen2},  there exists $H>0$ such that    $F(x,v) - Av^{\widetilde r}\geq0$ for $v>H$, $x\in B $, while   \pref{Hp_PS_cresc} implies that  $F(x,v) - Av^{\widetilde r}\geq -\widetilde C_0(H+H^{\widetilde \p}/{\widetilde \p})- AH^{\widetilde r} $ if $v\in[0,H]$. 
\\Then $F(x,v) \geq Av^{\widetilde r} -\widetilde C_0(H+H^{\widetilde \p}/{\widetilde \p})- AH^{\widetilde r} $ for $v\geq0$, $x\in B $.
\end{proof}

\subsection {Problems where $M$ is not a pure power}
We consider now a problem with the same right hand side of problem \pref{eq_rpq_p},  but where $M$ 
is not a pure power any more, which means that it interacts with the nonlinearity  in different ways near zero and at infinity. 
We prove the following result.
\begin{theorem}\label{th_notpure}
Let   $1<q<\p<p^*$ and 
suppose  $M$ satisfies \pref{Hp_PS_Mpos} and \pref{Hp_PS_SQii} with $\widetilde r<\p$.
Then the problem 
\begin{equation}\label{eq_rpq_M}
\left\{
\begin{array}{ll}
-M(\n{u}_W^p)\Delta_p u= -|u|^{q-2}u+\lambda|u|^{\p-2}u & \mbox{ in }\ \Omega\,,
\\
u=0 & \mbox{ on }\ \partial \Omega\,,
\end{array}
\right.
\end{equation}
has  at least one nonnegative nontrivial  solution    for  $\lambda>0$ small enough. Moreover,
\begin{enumerate}
\item   if  \pref{Hp_min_Morig} holds with  $r<p^*$, then the nonnegative nontrivial  solution exists for every $\lambda>0$,
\item  if  \pref{Hp_nomin_Mle} holds with $r>p^*$, then a further  nonnegative nontrivial solution exists for  $\lambda>0$ small enough.
\end{enumerate}
\end{theorem}
A model example for the nonlocal term $M$ could be 
$$M(s^p)=min\pg{s^{r_1-p},s^{r_2-p}},\qquad s\geq0,$$ with $r_2\in(0,\p)$, which guarantees \pref{Hp_PS_SQii}. 
Case 1. and case 2. in Theorem \tref{th_notpure} would correspond,  respectively, to  $r_1<p^*$ and to $r_1>p^*$.
\begin{proof}
We will find critical points of the functional $J^+$ in \pref{eq_funct_gen+} with $F(x,v)=-\frac1q|v|^q+\frac \lambda \p|v|^\p$.
The PS condition is satisfied by Proposition \tref{lm_PSgen}, since (\ref{Hp_PS_cresc}\,-\,\ref{Hp_PS_SQ}\,-\,\ref{Hp_PS_Mpos}) hold true.
Moreover, we can use the estimates in Lemma \tref{lm_estMF}.
Since 
 $J^+(u)\geq\frac1p \widehat M(\n u_W^p)-\lambda C \n{u}_W^\p$, by \pref{eq_est_ap3} there exist  $\Lambda,S,\rho>0$ such that 
\begin{equation}\label{eq_est_ap1}
J^+(u)\geq S\qquad\text{ for $\n{u}_W=\rho$ and $\lambda\in[0,\Lambda)$}.
\end{equation}
Having fixed $\lambda\in(0,\Lambda)$, let $\phi\in\W$ be a  positive function with $\n\phi_W=1$,  then  $J^+(T\phi)<0$ for some $T>\rho$ large enough, since $\lambda>0$, $\widetilde r<\p$ and then by \pref{eq_est_ap2}  the behavior of $J^+(t\phi)$ for large $t$ is given by the term with the exponent $\p$.

It is then possible, for $\lambda\in(0,\Lambda)$, to apply the Mountain Pass theorem to obtain a  nontrivial nonnegative solution, at level at least $S$.

When \pref{Hp_min_Morig} holds with $r<p^*$, 
 we can argue as for Corollary  \tref{coro_simin_model} to conclude that  the origin is a local minimum, this means that we can apply the Mountain Pass theorem even without the need for estimate \pref{eq_est_ap1},  then we obtain a nonnegative nontrivial solution for every $\lambda>0$ (actually, this case is also a consequence of Theorem \tref{th_notcoe_gen}). 

On the other hand, when   \pref{Hp_nomin_Mle} holds with $r>p^*$,  we are in the conditions of Theorem \tref{th_nomin_gen} and then  $J^+(w)<0$ for some  $w\in\W$ small in norm. 
As a consequence, for $\lambda\in(0,\Lambda)$, there exists a further nontrivial nonnegative solution which is a local minimum at a negative level in the ball $\pg{\n{u}_W\leq\rho}$. 
\end{proof}

\subsection{Some remarks about the applications.}\label{sec_commap}
\begin{itemize}
\item As mentioned before,  in some cases it is possible to guarantee that the nonnegative nontrivial solutions encountered are in fact strictly positive, using the maximum principle in \cite[Theorem 5]{Vazquez}, which we can apply to the local problem \pref{eq_main_gam}, also satisfied by the solutions.  This is the case in the theorems \tref{th_coe_mod}, \tref{th_notcoe_mod}, \tref{th_notcoe_mod_mult} and \tref{th_notpure}, if we take $q\geq p$, while in the theorems \tref{th_coe_gen} and \tref{th_notcoe_gen} one needs to assume $f(x,v)\geq -v^{q-1}$ for $v$ positive and small, with $q\geq p$. 

\item
It is  interesting to observe the behavior of the  problem \pref{eq_rpq_p} with   $q<p<\p$ in function of the nonlocal term, controlled by the parameter  $r$:  in the local case $r=p$ there exists  a solutions for every $\lambda>0$ (see Theorem \tref{th_notcoe_mod} with $\mu=0$), however, as $r$ becomes larger than $\p$, existence is lost for $\lambda$ large, but it is recovered if $r>p^*$ (see Theorem \tref{th_coe_mod}).

The reverse happens if we consider  problem \pref{eq_rpq_p} with   $q<\p<p$:  in the local case $r=p$ there exist no solution for $\lambda$ large  (see Theorem \tref{th_coe_mod}), but when  $q<r<\p$  a solution exists for every $\lambda>0$ by Theorem \tref{th_notcoe_mod}.

\item
If one considers the local problem  $r=p$ in Theorem \tref{th_coe_mod}, since $\p<r$,  the nonlinearity is $p$-sublinear.  Such problem was considered in \cite{Ane12_2nl_blw} and  the availability of the sub and supersolutions method allowed the author to obtain the analogous of the first point of Theorem \tref{th_coe_mod}, but with  a better description of the set of parameters for which either two, one, or zero solutions exist. 
 Our result for the nonlocal problem turns out to be less precise in view of the fact that the sup and supersolutions method cannot be used as in the local case (see \cite{ItuGM_contra,FiSu18_Kirch_CP_contra,itumas_contrex}).

On the other hand, our result shows that for the nonlocal problem, a new behavior arises  as $r$ goes over the threshold $p^*$, even if the interaction with the nonlinearity seems unchanged, since  the geometry of the functional changes near the origin, which is not a local minimum any more. The same phenomenon can be see in Theorem \tref{th_notpure}.

\item 
If  $r=p$ in  problem \pref{eq_rpqs}, one obtains, after a rescaling to match the different parametrization used there, the problem considered in \cite{Ane10_3pow}, where a  three solutions result similar to our Theorem \tref{th_notcoe_mod_mult} is obtained. 
Also, Theorem \tref{th_notcoe_mod} in the case   $r=p=\s=2$ was already proved in  \cite{MIO_ODA}.

\item
Our  results are also related with those in \cite{AmbArc_Kirch_Bif,AmbArc_Kirch_VarDeg}, where problems similar to \pref{eq_rpq_M} were considered for $p=2$ and with various  hypotheses on the nonlocal term $M$, including the possibility of being degenerate. In particular, the nonlinearity  did not include the term $-|u|^{q-2}u$, but it was also considered an additional  forcing term $h\in L^2$. 
\end{itemize}

\section{A remark about a-priori estimates}\label{sec_CE_Lieb}
In this section we present an  example that shows that the nonlocal term may make it impossible to obtain certain uniform estimates that hold true for the local version of equation \pref{eq_main_W_f}.

In \cite{AlvCor_Kirch2} it was observed that, in the nondegenerate case where $M\geq m_0>0$, standard regularity results may be used since, as we noted in Section \tref{sec_NotRem},  if $u$ is a solution of \pref{eq_main_W_f} then it is also  a solution of problem \pref{eq_main_gam} with $\llambda=1/M(\n{u}_W^p)$, and in this case $\llambda$ is bounded by $1/m_0$. In the degenerate case it is still possible to use this argument in order to guarantee, using for instance the results in  \cite{Guedda_Veron},  that  solutions $u$ are always of class  $\Cont^{1,\alpha}$, for some $\alpha>0$, provided $M(\n{u}^p_W)\neq0$. 
However,  it is not clear if one can obtain a uniform estimate for the $\Cont^{1,\alpha}$ norm of the solutions, as in \cite{Lieberman}, due to the fact that the multiplier $\llambda$ can be arbitrarily large.   Because of the lack of this uniform estimate,  in  Theorem \tref{th_simin}, we had to assume our minimum to be with respect to the $L^\infty$ norm, while in Theorem \tref{th_simin_C1}, for the stronger result that only assumes the minimum to be in the $\Cont^1$ norm, we needed the more restrictive  assumptions (\ref{eq_condr_C1H1}\,-\,\ref{eq_condlt_C1H1}), which in fact give a restriction on how fast the multiplier $\llambda=1/M(\n{u}_W^p)$ in \pref{eq_main_gam} can become large, in comparison with the behavior of the nonlinearity.

Below we actually show a situation where there is not such a bound as the  one that is obtained, in the local case,  by the results in \cite{Lieberman}. The example considers the Laplacian case ($p=2$) and a power behavior for the nonlocal term $M$, but it shows that one should not expect, also in more  general settings, to be able to obtain such uniform bounds in the presence of a degenerate term of Kirchhoff type.
In fact, for $p=2$ one would even expect to obtain better regularity results than for the quasilinear case.

\par\medskip \newcommand{\M}{D}

Consider the nonlocal problem  \pref{eq_PMq}
\begin{equation}\tag{$P_a$} \label{eq_PMq}
\begin{cases}
-M(\n{u}_W^2)\Delta u=g_a(u)=-au^{q-1}+u^{\p-1}&in \ \Omega\,,\\
u=0&on \ \partial \Omega\,,
\end{cases}
\end{equation}
 with parameter $a\in(0,A]$ and suitable  $1<q<\p<2$, where again $\Omega\subset\mathbb{R}^{N}$ is a bounded and smooth domain.
In the case where  $M$ is constant it is possible to  apply the results in 
\cite{Lieberman}: in fact, given $\M,A>0 $, for a suitable constant $\Lambda$, the right hand side satisfies   $|g_a(s)|\leq \Lambda$ for every $s\in [-\M,\M]$, $a\in(0,A]$. Then by \cite[Theorem 1]{Lieberman}  there exist
 $\beta\in(0,1)$ depending on $\Lambda,N$, and $C>0$ depending on $\Lambda,\M,N,\Omega$, such that$$\n{u}_{\Cont^{1,\beta}}\leq C\,$$ for any weak solution satisfying $\n{u}_\infty<\M$.

The result in the following proposition  shows that this  uniform estimate is false for every $\beta\in(0,1)$, and even in $\Cont^1$,  for  suitable degenerate nonlocal terms  $M$.
\begin{prop}\label{prop_priori}
If $M(s^2)={s^{r -2}}$
with $r \in(2+\frac2{N },2^*)$, $N\geq3$, 
then there exists  a family of functions, satisfying  problem \pref{eq_PMq} with  $ a\in(0,1]$, which is bounded in $L^\infty$ but unbounded in $\Cont^1$.
\\
In the limiting case $r =2^*$, a family as above exists, satisfying  problem \pref{eq_PMq} with a fixed  $ a>0$. 
\end{prop}
\begin{proof}
Let $B$ be an open ball compactly contained in $\Omega$. Without loss of generality we suppose $B$ is centered at the origin and has radius $1$.

In \cite{IlEg10_Hopfviol} it was proved that for  $\p>q>1$ small enough and  $b_0>0$ large enough, there exists a nonnegative solution $\Phi\in\Cont^{1,\beta}$, $\beta>0$,   of the problem 
$$\begin{cases}
-\Delta u=- u^{q-1}+b_0u^{\p-1}&in \ B,\\
u=0&on \ \partial B,
\end{cases}$$
which has compact support, and then can be continued by zero to a solution on the whole of $\Omega$. 
We consider  now, for $\lambda\geq 1$, $\mu\in(0,1]$, the family of functions defined in $\Omega$ as $$\Phi_{\lambda,\mu}(x):=\begin{cases}
\mu\Phi(\lambda x) & in\ B_{1/\lambda}\,,\\
0& in\ \Omega\setminus B_{1/\lambda}\,,
\end{cases}$$
where $B_{1/\lambda}$ is  the ball centered at the origin with radius $1/\lambda$.
Then $\Phi_{\lambda,\mu}$  satisfies the problem 
$$\begin{cases}
\displaystyle -\frac{\Delta u}{\mu^{2-\p}\lambda^2b_0}={-\frac{\mu^{\p-q}}{b_0}u^{q-1} + u^{\p-1}}&in \ \Omega\,,\\
u=0&on \ \partial \Omega\,.
\end{cases}$$
Straightforward computations show that 
\begin{equation}\label{eq_CE_2} \begin{cases} 
\n{\Phi_{\lambda,\mu}}_\infty&=\mu \n{\Phi_{}}_\infty\,,
\\\n{\nabla\Phi_{\lambda,\mu}}_\infty&=\mu\lambda \n{\nabla\Phi_{}}_\infty\,,
\\\n{\Phi_{\lambda,\mu}}_W^2&=\mu^2\lambda^{2-N}\n{\Phi_{}}_W^2\,,
\end{cases}
\end{equation}
moreover the functions $\Phi_{\lambda,\mu}$ are solutions of the nonlocal problem \pref{eq_PMq} provided 
\begin{equation} \label{eq_CE_1} \begin{cases} 
M(\n{\Phi_{\lambda,\mu}}_W^2)={\pt{\mu^2\lambda^{2-N}\n{\Phi_{}}_W^2}^{(r -2)/2}}=(\mu^{2-\p}\lambda^2b_0)^{-1}\,,\\
  a=\mu^{\p-q}/b_0\,.
  \end{cases}
  \end{equation}
If we take  $\mu(\lambda)=\lambda^{-\alpha }/E$ with $\alpha \in[0,1)$, $E>0$, then \pref{eq_CE_2} shows that the family $\pg{\Phi_{\lambda,\mu(\lambda)}}_{\lambda\geq1}$ is bounded in $L^\infty$  but   unbounded in $\Cont^1$ and \pref{eq_CE_1} becomes
\begin{equation}\label{eq_issolLie}\begin{cases}
{\pt{\lambda^{2(1-\alpha )-N}}^{(r -2)/2}}\lambda^{2+\alpha (\p-2)}(b_0\n{\Phi_{}}_W^{r -2}/E^{r -\p})=1\,,
\\a=\pt{\frac{\lambda^{-\alpha }}{E}}^{\p-q}/b_0\,.
\end{cases}
\end{equation}
From \pref{eq_issolLie} we see that we can choose
$$E=\pt{\n{\Phi_{}}_W^{r -2} b_0}^{1/(r-\p) }\qquad \text{and}\qquad
 (r -2)(1-\alpha -N/2)=-2+\alpha (2-\p)\,, $$
 that is 
$$\alpha=\frac{N+r-rN/2}{r-\p}=\frac{N}{2^*}\frac{2^*-r}{r-\p}\,.$$
Since we need $\alpha\in[0,1)$, we obtain $\p<r\leq 2^*$  and $r>2^*\frac{N+\p}{N+2^*}=2\pt{1+\frac\p N}$.  Remembering that in \cite{IlEg10_Hopfviol}  $\p>1$ could be taken as near to $1$ as wanted, we get the condition in the claim that $2+\frac{2}{N}<r\leq 2^*$.

 Observe that $a<1$ for $\lambda$ large, except for the case $r =2^*$ where $a$ is constant. We have thus proved our claim.
\end{proof}

\begin{remark}
It is interesting to compare the above proof with the conditions assumed in Theorem \tref{th_simin_C1}, see also Remark \tref{rmk_l_lt}. Actually, here we had to take $\p-1$ and $q-1$ very small, which goes in the opposite direction with respect to the   assumptions (\ref{eq_condr_C1H1}\,-\,\ref{eq_condlt_C1H1}).
\eor
\end{remark}
\begin{remark}
Other examples, involving nonlocal terms that are degenerate at infinity, can be easily obtained.
For instance, 
let $\Omega=(0,\pi)$ and $\phi_i=\sqrt{\frac2\pi}\sin(ix)$ be the normalized eigenfunctions of the Laplacian associated to the eigenvalue $i^2$:
$$
\left\{
\begin{array}{lcl}
- \phi_i'' =  i^2 \phi_i& {\rm in} & (0,\pi)\,,\\
\phi_i (0)=\phi_i(\pi)= 0\,.  & \;
\end{array}\right.
$$
Then they are a family of functions which is bounded in $L^\infty$ and unbounded in $\Cont^1$, which are solutions of the nonlocal problem, with weight $M(s^2)=s^{-2}$,
$$
\left\{
\begin{array}{lcl}
- \frac1{\n{u}^2_W} u'' =   u& {\rm in} & (0,\pi)\,,\\
u (0)=u(\pi)= 0\,.  & \;
\end{array}\right.
$$
\eor
\end{remark}

\section{Appendix}
In this Appendix we summarize some estimates, which were used in the proof of Theorem \tref{th_nomin_gen}.

As is classical in the literature, we consider  compact support approximations of the instanton functions (see \cite{Azor_Alo})
\begin{equation}\label{eq_def_phi}
\Phi_{\epsilon}^{}(x)=\left(C_{N,p}\ \frac{\epsilon^{\frac1{p-1}}}{\epsilon^{\frac p{p-1}}+|x|^{\frac p{p-1}}}\right)^{\frac{N-p}p}\ \ \ \ \ \ \ \ \ \ \mbox{with $\epsilon >0$},
\end{equation}
which realize the best  constant $S$ for the  Sobolev embedding inequality $ S\n{\Phi_{\epsilon}^{}}_{L^{p^*}(\R^N)}^{p}\leq  \n{\nabla U}_{L^{p}(\R^N)}^p \,,$
where the constant $C_{N,p}$ is  chosen in such a way that $\n{\nabla\Phi_{\epsilon}^{}}_{L^{p}(\R^N)}^p=\n{\Phi_{\epsilon}^{}}_{L^{p^*}(\R^N)}^{p^*}=S^{N/p}$.

In order to obtain a better behavior of the estimates of certain norms, we use here a technique initially proposed in \cite{GazRuf}, see also \cite{CalancRuf}: it consists in using a further parameter that controls the size of the support of the cutoff function: 
we take $\xi_{m}\in C^{\infty}_{0}(\mathbb{R}^{N})$ such that $0\leq\xi_{m} (x)\leq 1$, $||\nabla\xi_{m}||_{\infty}\leq 4m$ and 
\begin{equation*}
\xi_{m} (x)=
\begin{cases}
1, & \mbox{if $x\in B_{\frac{1}{2m}}$,}\\
0, & \mbox{if $x\in\mathbb{R}^{N}\setminus B_{\frac{1}{m}}$,}
\end{cases} 
\end{equation*}
where $B_{r}$ is  the ball centered at the origin with radius $r$.
We suppose $0\in\Omega$ and we 
define $\Phi_{\epsilon,m}(x)=\xi_{m} (x)\Phi_{\epsilon}^{}(x)|_\Omega$, so that for $m$ large enough, $\Phi_{\epsilon,m} \in \W(\Omega)$. 
Then one obtains the  estimates contained in the following lemma.
\begin{lemma}\label{lm_GaRuCa_p}
Suppose that  $m\rightarrow\infty$ (or $m$ is constant but large enough),  $\epsilon\rightarrow 0$ and $\epsilon m\to 0$. Then the following estimates hold for suitable constants $D_1,..,D_4>0$ having the indicated dependencies.
\begin{enumerate}[label=(\alph*)]
  \item \label{p1_p} $\left| ||\Phi_{\epsilon,m}||^{p}_{W}-S^{\frac{N}{p}}\right|\leq D_{1}(N,p)(\epsilon m)^\frac{N-p}{p-1}$
  \item \label{p2_p} $\left| ||\Phi_{\epsilon,m}||^{p^{*}}_{p^{*}}-S^{\frac{N}{p}}\right|\leq D_{2}(N,p)(\epsilon m)^\frac{N}{p-1}$
\item  \label{p36_p} 
$D_3(N,p,s) \eta_{N,p,s}(\varepsilon,m)\leq ||\Phi_{\epsilon,m}||^{s}_{s}  \leq D_4(N,p,s)\eta_{N,p,s}(\varepsilon,m)$
\\where 
$$\eta_{N,p,s}(\varepsilon,m)=\begin{cases}
                \epsilon^{N-\frac{N-p}{p}s},                 &\ \ \ \mbox{if $p^*>s>\frac{N(p-1)}{N-p}$}, \\
\epsilon^{N-\frac{N-p}{p}s}|\log(\epsilon m)|{=\epsilon^{N/p}|\log(\epsilon m)|,   }&\ \ \ \mbox{if $s=\frac{N(p-1)}{N-p}$},\\
 \epsilon^{N-\frac{N-p}{p}s}\pt{\epsilon m}^{s\frac{N-p}{p-1}-N}  {=\epsilon^{\frac{N-p}{p(p-1)}s} m^{s\frac{N-p}{p-1}-N}}                                                      &\ \ \ \mbox{if $1\leq s<\frac{N(p-1)}{N-p}$}.                                
\end{cases}$$
\end{enumerate}
\end{lemma}

The additional parameter $m$ can be used in order  to improve some of the estimates. We will set $m=\varepsilon^{-\beta}$ with $\beta\in(0,1)$ suitably near to $1$ and define 
\begin{equation}\label{eq_defpsi_p}
\psi_\varepsilon=\Phi_{\epsilon,\varepsilon^{-\beta}}.
\end{equation}
Then  one easily obtains the following estimates.
\begin{corollary}\label{cor_GaRuCa_p} 
For any  $\ttau>0$, there exists $\beta\in(0,1)$ (depending on $N, p, s$) such that, when $\varepsilon\to0$,
\begin{equation}\label{eq_estPsi_H_p}
\left| ||\psi_\varepsilon||^{p}_{W}-S^{\frac{N}{p}}\right|\leq D_{1}(N,p)(\epsilon^{1-\beta})^\frac{N-p}{p-1}\to0\,,
\end{equation}
  \begin{equation}\label{eq_estPsi_p*_p} \left| ||\psi_\varepsilon||^{p^{*}}_{p^{*}}-S^{\frac{N}{p}}\right|\leq D_{2}(N,p)(\epsilon^{1-\beta})^\frac{N}{p-1}\to0\,,
\end{equation}
\begin{align}\label{eq_estPsi_s_p}
 D_{3}(N,p,s)\epsilon^{N-\frac{N-p}{p}s}\leq||\psi_\varepsilon||^{s}_{s} & \leq D_4(N,p,s)
\begin{cases}
\epsilon^{N-\frac{N-p}{p}s}                                                                    &\ \ \ \mbox{for $p^*>s>\frac{N(p-1)}{N-p}$}, \\
\epsilon^{N-\frac{N-p}{p}s-\ttau}      &\ \ \ \mbox{for $1\leq s\leq\frac{N(p-1)}{N-p}$}.
\end{cases}
\end{align}
\end{corollary}
\begin{remark}
The proof of the estimates in Lemma \tref{lm_GaRuCa_p} can be obtained following the lines of the classical arguments used in \cite{BrezNir_crit}, see also \cite{Ghous_book}.

The estimates \ref{p1_p} and \ref{p2_p} appear in \cite[Lemma 6]{GazRuf} for the case $p=2$ and can be easily generalized to the general case $p>1$, where, for constant $m$, they reduce to those in  \cite[page 947]{Azor_Alo}. 

The lower estimate in \ref{p36_p}, for a constant $m$, is analogous to the one obtained in \cite[Lemma A5]{Azor_Alo}, while the proof of the upper estimate in \ref{p36_p} can be seen in \cite[Lema 5.5]{teseM_Mayc}  for the case $p=2$, and we briefly sketch it below for sake of completeness. 
In fact, when  $m\to\infty$, a correcting term in $m$ appears when $(\Phi_{\epsilon})^{s}$ fails to be  integrable in $\R^N$, that is, when $-\frac{N-p}{p-1}s+N-1\geq-1$, because the reduction in the support of $\Phi_{\epsilon,m}$ reduces the contribution of $\Phi_{\epsilon}^{s}$ away from the origin. \eor
\end{remark}

\begin{proof}[Sketch of the proof of point \ref{p36_p} in Lemma \tref{lm_GaRuCa_p}]
For $m$ large enough one has
\begin{align}\label{eq_phisp_1}
    \int_{B_{\frac{1}{2m}}}|\Phi^{}_{\epsilon}(x)|^{s}dx           \leq &||\Phi_{\epsilon,m}||^{s}_{s}    \leq \int_{B_{\frac{1}{m}}}|\Phi^{}_{\epsilon}(x)|^{s}dx\,.
\end{align}
By the change of variable  $x=\epsilon w$, for $k=1,2$, we get
\begin{align}\label{eq_phisp_2} \int_{B_{\frac{1}{km}}}\left(\frac{\epsilon^{\frac1{p-1}}}{\epsilon^{\frac p{p-1}}+|x|^{\frac p{p-1}}}\right)^{\frac{N-p}ps}dx=\epsilon^{N-\frac{N-p}{p}s}\int_{B_{\frac{1}{k\epsilon m}}}\frac{1}{(1+|w|^{\frac p{p-1}})^{\frac{N-p}{p}s}}dw,
\end{align}
where we already obtain the  term $\epsilon^{N-\frac{N-p}{p}s}$, common to all the estimates in point \ref{p36_p}.
We need to estimate the  last integral in \pref{eq_phisp_2}. If it converges, as $\epsilon m\to 0$, then it can be bounded between two constants.
If it does not converge, then an additional factor appears, whose asymptotic behavior is 
\[\begin{cases}
               |\log(\epsilon m)|&\ \ \ \mbox{if $s=\frac{N(p-1)}{N-p}$},\\
\pt{\epsilon m}^{s\frac{N-p}{p-1}-N}                                         &\ \ \ \mbox{if $1\leq s<\frac{N(p-1)}{N-p}$}.                         
\end{cases}\qedhere\]
\end{proof}

\section*{Acknowledgement}

\noindent Funding: 

\noindent
L. Iturriaga was partially supported by Programa Basal PFB 03, CMM, U. de Chile;
Fondecyt grants 1161635, 1171691 and 1181125 (Chile).

%\vspace{.3cm}

\noindent E. Massa was  supported by: grant $\#$2014/25398-0, São Paulo Research Foundation (FAPESP) and  grants $\#$308354/2014-1 and $\#$303447/2017-6, CNPq/Brazil.

\providecommand{\bysame}{\leavevmode\hbox to3em{\hrulefill}\thinspace}

\end{document}